\documentclass[12pt,twoside,a4paper,notitlepage]{article}

\usepackage{graphicx}
\usepackage{amsthm}
\usepackage{mathtools}
\usepackage{enumerate}
\usepackage{amssymb}
\usepackage{mathptmx}     
\usepackage{csquotes}

\usepackage[square, authoryear]{natbib}
\bibpunct{[}{]}{;}{t}{,}{,}
\setcitestyle{square,aysep={},yysep={;},citesep{;},notesep={, }}
\usepackage{latexsym}

\begin{document}


\title{Bernard Bolzano: from Topological to Arithmetical Continuum and Back Again}

\author{Kate\v{r}ina Trlifajov\'{a}}

\date{}







\maketitle

\begin{abstract}

Although Bolzano's concept of the continuum has gradually evolved, the basis remained the same: the continuum as an infinite class of points arranged in such a way that the so-called \emph{Bolzano completeness} holds. Bolzano realized over time that  the central role to a general comprehension of continuum plays its arithmetic description and constructed his measurable numbers. Their interpretations in the standard and non-standard models of real numbers clarify their relationship and also suggest why Bolzano did not base his theory of functions on infinitesimal numbers. The three main theorems on measurable numbers are actually various forms of their completeness. I argue why the second one is indeed the \emph{Supremum Theorem} and that an important corollary of the third one is a proof of the \emph{Bolzano completeness}. Only when the notion of continuum was supported by measurable numbers could Bolzano, in his last book, \emph{Paradoxes of the Infinite}, confidently defend the general properties of the continuum and reject the paradoxes associated with them.


\end{abstract}
\section{Preface}

Bolzano was concerned with the question of the continuum throughout all his life, from his first mathematical text \emph{Considerations on Some Objects of Elementary Geometry} in 1804 
to his last \emph{Paradoxes of the Infinite} in 1848.
From the beginning he conceived of the continuum set-theoretically and topologically; as an infinite collection of points arranged in a particular way. This conception has remained essentially unchanged, only being refined and supplemented. Bolzano gradually comprehended that a full understanding of continua required an arithmetic description of them, which did not exist at the time; both Dedekind and Cantor theory of real numbers were not developed until 1872. Bolzano constructed his measurable numbers, about which all properties of real numbers hold, in the late 1830s and proved that they satisfy the arrangement of his topological continuum. 

\subsection{The content} 

The content of this article is as follows. In Section \ref{long}, I summarize the  evolution of Bolzano's topological concept of the continuum over time while referring to the works of Bolzano's scholars. 
Measurable numbers that represent Bolzano's arithmetization of continuum are briefly introduced in Section \ref{measurable}. Their interpretation in Section \ref{interpretation} shows the relation of measurable numbers to the standard and non-standard model of real numbers once called for by \citet{lakatos1978Cauchy}. Section \ref{main} deals with three main Bolzano's theorems on measurable numbers and demonstrates that they correspond to four forms of completeness of real numbers. The first theorem generally accepted as \emph{Bolzano-Cauchy's} is left aside but the second one, variously called the \emph{Upper-Bound}, the \emph{Least-Upper-Bound}, or the \emph{Greatest-Lower Bound}, is discussed, as well as the third theorem, which seems to be somehow related to \emph{Dedekind}'s cuts. I believe that the main purpose of the last one was to demonstrate that the ordering of measurable numbers is that of the continuum. Finally, in the last Section \ref{back} I return  to the Bolzano's mature concept of continuum to show the influence of measurable numbers on it.

\subsection{Abbreviations}

I quote Bolzano's work according to the translation contained in \emph{The Mathematical Work of Bernard Bolzano} \citep{russ2004Mathematical}  and use common abbreviations, derived from German titles of works, and numbers of paragraphs.


\begin{itemize} 
\item[BG] \emph{Considerations on Some Objects of Elementary Geometry}, [pp. 25 - 82].

\emph{Betrachtungen über einige Gegenstände der Elementargeometrie}, Prague 1804. 

\item[BD] \emph{Contributions to a Better-Grounded Presentation of Mathematics}, [pp. 83 - 138].

\emph{Beytr\H age zu einer begr\H undeteren Darstellung der Mathematik}, Prague 1810.

\item [RB] \emph{Purely Analytic Proof of the Theorem, that between any two Values which
give Results of Opposite Sign, there lies at least one real Root of the Equation}, [pp. 251 - 278].

\emph{Rein analytischer Beweis des Lehrsatzes, dass zwischen je zwey Werthen,
die ein entgegengesetztes Resultat gew\H ahren, wenigstens eine reelleWurzel
der Gleichung liege}, Prague 1817.

\item [DP] \emph{The Three Problems of Rectification, Complanation and Cubature solved
without consideration of the infinitely small, without the hypotheses of
Archimedes, and without any assumption which is not strictly provable}, [pp. 279 - 344].

\emph{Die drey Probleme der Rectification, der Complanation und der Cubirung
ohne Betrachtung des unendlich Kleinen, ohne die Annahmen des
Archimedes, und ohne irgend eine nicht streng erweisliche Voraussetzung
gel\H ost.  zugleich als Probe einer g\H anzlichen Umstaltung der Raumwissenschaft,
allen Mathematikern zur Pr\H ufung vorgelegt}, Leipzig 1817.

\item[RZ] \emph{Pure Theory of Numbers. Seventh Section: Infinite Quantity Concepts}, [pp. 355 - 428]. 

\emph{Reine Zahlenlehre
Siebenter Abschnitt. Unendliche Größenbegriffe}, BGA 2A8 Stuttgart, 1976.

\item[PU] \emph{Paradoxes of the Infinite}, [pp. 591 - 678]. 

\emph{Paradoxien des Unendlichen}, Leipzig 1851. 

\end{itemize}

\section{The topological concept}\label{long}

The concept of continuity remained a subject of investigation throughout Bolzano's entire scientific career.  

\begin{quote} Bolzano's interest in his fundamental geometrical problem of defining the concepts of line, surface, solid, and continuum spanned his entire creative lifetime, from his earliest thoughts on mathematics and first published works to his last research papers. \citep[p. 294]{johnsondalem.PreludeDimensionTheory1977}. 
\end{quote}
 
\subsection{The early work}

Already in his first published mathematical work, \emph{Considerations on Some Objects of Elementary Geometry} [BG] from 1804, he tries to formulate a \enquote{proper theory of the straight line} - its existence, its determination, its possible infinite extension - as a part of \emph{pure} geometry. All geometrical doctrines concerning triangles, parallels, etc. should be derived from this basis. 

Bolzano opens his next published 1810 work, \emph{Contributions to a Better-Grounded Presentation of Mathematics} [BD] by the assertion that although mathematics
\enquote{comes nearest to the ideal of perfection \dots some gaps and imperfection are still to be found even in the most elementary theories of all mathematical disciplines.} (BD Preface). This particularly applies to geometry, where precise definitions of line, surface and solid are still lacking. Bolzano builds these concepts on a conviction he has held all his life. 
\begin{quote} The straight line, the plane are objects of a composite kind, in which, we imagine, for instance, innumerably many points as well as particular relationship which these points must have to some given ones. [BD II. \S 5]. 
\end{quote} 

Rusnock points out that this concept of the continuum was influenced by Baumgarten's \emph{Metaphysik} from 1783, the first philosophical book Bolzano had read at the age of sixteen. Baumgarten claims that extensions belong to composite things and only extended things have shape. Monads as simple beings have neither shape nor extension, though wholes composed of them can have both. So too in geometry, where multitudes of unextended points can constitute extended wholes. Though Bolzano rejected some errors he identified in Baumgarten's \emph{Metaphysik} - for instance that continuum could be composed of only finitely many points - he accepted and defended these ideas. \citep[p. 194]{shapiroHistoryContinuaPhilosophical2021}. 

In a polemic regarding the distinction between empirical and a priori judgements Bolzano rejects Kantian doctrine of constructing objects through \emph{pure intuition} [BD I. \S 6] and announces a program for grounding mathematics on a purely analytic foundation. See  \citep[pp. 45 - 50]{rusnock2000Bolzanos}.

Later, Bolzano developed his own concept of physical continuum, a \enquote{monadology} in which he proposed the existence of infinitely many partless substances, which he calls \emph{atoms}, completely filling space, and interacting via forces of attraction and repulsion. 
\citep{simons2015Bolzanos}.

\subsection{Geometrical continuum}
In his 1817 work \emph{Three problems of the Rectification, Complanation and Cubature} [DP],  Bolzano pursues a rigorous mathematical definition of a \emph{line} that should serve as the foundation for the entire structure of analysis. He found all other previous definitions inadequate, particularly those that relied on infinitesimals or the Archimedean axiom.

Bolzano's definition is based on a set-theoretic comprehension. In [DP \S 11], he defines a \emph{spatial object} as a system (a collection) of points. In addition to the problem of continuity, he also deals here with the question of dimension. A line arises as a \emph{composition of points}, it requires an infinite multitude of points which are arranged in a special way:

\begin{quote} A spatial object, at every point of which, beginning at a certain
distance and for all smaller distances, there is at least \emph{one} and at most only a \emph{finite} set of points as neighbours, is called a \emph{line in general}. 
\end{quote}

The meaning of the terms is better clarified from his comments and from diagrams of examples he gives.  
\begin{enumerate}
\item A \emph{distance} is a \enquote{determinate kind of separation of two points}. [BG \S 24]. It can be determined in any way, geometrically or algebraically.
\item A \emph{neighbour} in a given distance is a point of the intersection of a spatial object with the surface of a sphere with radius equal to the given distance 
\citep[p. 276]{johnsondalem.PreludeDimensionTheory1977}. 
\item \emph{Lines} can be curve lines that can intersect. Hence there can be more neighbours at the same distance.

\item At most a \emph{finite} set of points, otherwise it might be a point of a surface or of a solid.
\end{enumerate}

\begin{figure}[h]
\centering
\includegraphics[scale=0.75]{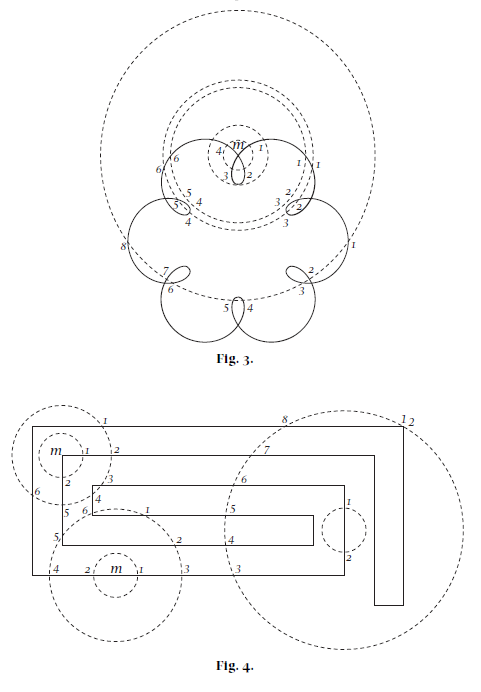}
\caption{Bolzano's picture}\label{fig1}
\end{figure}

Bolzano continues and defines with help of the number of neighbours in a given distance such notions as \emph{boundary} and \emph{internal} points, a \emph{bounded}, \emph{simple}, \emph{closed} line, and \emph{connected} line.  A \emph{surface} and a \emph{solid} is defined on the same principle, .  

\begin{quote} For with the surface and the solid this is different: in the former at each point there is a complete line of points and in the latter even a complete surface of points which have the same distance from the given point. [DP \S 12.] 
\end{quote}

Bolzano maintained this concept surprisingly modern for his time throughout whole his life, although he slightly modified its verbal terms. Nearly the same definitions occur in the [PU \S 49]. 

\begin{quote} The topological and set-theoretical methods which Bolzano used as a way of dealing with the basic elements of geometry were without question far ahead of his time. His topological basis, derived from his concept of neighbour and then later also from his concept of isolated point, is very deep. \citep[p. 292]{johnsondalem.PreludeDimensionTheory1977}.  \end{quote}

Compared to Cantor’s later definition, Bolzano’s concept of a continuum does not necessarily imply connectedness. Two separate lines can be continuous yet not connected.

\subsection{Arithmetical continuum}\label{continuous}

That same year 1817, Bolzano published \emph{Purely Analytic Proof} [RB], in which he sought to provide an analytic, i.e. arithmetic, description of the intersection point of a line, the x-axis, and a curve, given by a continuous function, that crosses it. To prove it he gives his own \enquote{correct definition} of \emph{continuity of a function}. 

\begin{quote} The expression that a function $f(x)$ varies according to the law of continuity for all values $x$ inside or outside certain limits means only that, 
if $x$ is any such value the difference $f (x + \omega) - f(x)$ can be made smaller than any given quantity, provided $\omega$ can be taken as small as we please.\footnote{Bolzano explains he uses the phrase \emph{a quantity which can become smaller than any given quantity} or \emph{which can become as small as we please} instead of \emph{infinitely small quantities} with equal success in \emph{Preface} of BL \citep[p. 158]{russ2004Mathematical}}
 [\emph{Preface} RB], \citep[p. 256]{russ2004Mathematical}.\footnote{This famous Bolzano's definition of continuity has commonly been considered as an early presentation of the modern – Weierstrassian – epsilon-delta style. \citet{fuentesguillen2020Notion} point out subtleties of used notions. 
}\end{quote} 

This definition has an important consequence that Bolzano repeatedly used in the proof which we will return in Section \ref{IMV}.. 

\begin{quote}If functions $f(x), \varphi(x)$ vary according to the law of continuity then by virtue of the law of continuity $$(\ast) \quad f(\alpha) < \varphi(\alpha) \Rightarrow f(\alpha + i) < \varphi(\alpha + i)$$ 
provided $i$ is taken sufficiently small. \emph{Preface} of BL, \citep[p. 260]{russ2004Mathematical}.
\end{quote}

However, Bolzano encountered a conceptual obstacle at the time: the absence of an arithmetization of continuum, as no rigorous theory of real numbers yet existed. 
Despite this, he successfully formulated what is now known as the \emph{Cauchy} or rather the \emph{Bolzano-Cauchy} criterion for the convergence of a sequence (from now on BC-criterion).
\footnote{Bolzano states BC-criterion in [RB \S 7] four years before Cauchy in his \emph{Cours d'Analyse}.  
There has been a discussion on the possibility that Cauchy might have plagiarized from Bolzano. \citep[p. 149]{russ2004Mathematical}. However, \citep{bair2020Continuity} argue that the two concepts were developed independently.} 
He also stated and proved the \emph{Supremum Theorem} [RB \S 12], see Section \ref{supremum}. Then he could finally prove the \emph{Intermediate Value Theorem} in [RB \S 15], see Section \ref{IMV}.

\subsection{Later writings}

Two additional unfinished essays on the definition of extension and the continuum can be found among Bolzano’s manuscripts  \cite[p. 199 - 202]{rusnock2020Bolzano}. In \emph{Attempted Definitions of the Concepts of a Line, Surface, and Solid}, Bolzano

\begin{quote} calls an element $p$ of a spatial object \emph{connected} iff for some positive number d, there is at least one element of the spatial object (neighbour) lying at the distance h from p, for all $h \leq d$. A spatial object is extended, if all its points are connected.  \cite[p. 201]{rusnock2020Bolzano}. \end{quote} 

The idea of this definition is similar as in [DP]. However, this is more general, an extended spatial objects can be a line, a surface or a space. This property is called the \emph{Bolzano completeness}, see Section \ref{main}: 

\begin{quote}
Every element beginning at a certain distance and for all smaller distances has a neighbour. \end{quote}

Bolzano understands more and more that the central role for the comprehension of the general phenomenon of continuum plays the theory of quantities. The structure of the real (that he calls measurable) numbers so becomes the ultimate basis of the spatial, temporal, and the material continua.  \cite[p. 208]{rusnock2020Bolzano}.

 \begin{quote}We can see that it was metric intuition, interpreted in its narrowest sense, that guided his construction of the concept: distances are ultimately involved in his definition of the continuum in the form of real numbers. \citep[p. 522]{granger1969Concept}.\footnote{\enquote{On voit done que c'est l'intuition métrique, interpretée en son sens le plus étroit, qui a guidé sa construction du concept : les distances interviennent en fin de compte dans sa définition du continu sous la forme de nombres réels.}}
\end{quote}

Although Granger seems to regret Bolzano's abandonment of a purely topological definition of continuity, \citet[p. 208]{rusnock2020Bolzano} highlights that Bolzano sought with a fair measure of success to provide what was necessary: an account of the real numbers that is independent of the theories of space, time, and any of the other objects where the theory of real numbers can
be applied.

\section{Measurable numbers}\label{measurable}


Bolzano apparently attempted to remedy the lack of arithmetical descriptions of points of a line in the  VII\textsuperscript{th} of his extensive work \emph{Pure Theory of Numbers}, titled \emph{Infinite Quantity Concepts} [RZ], written in the 1830s.  Bolzano had never finished this work. It was preserved as a manuscript not ready for publication in Vienna's State Library until the second half of the last century \citep{rychlik1957Theorie}. The second Berg's critical edition including important supplements was published later \citep{bolzano1976Bernard}.



 The problem is that there occur some minor errors and  inconsistencies. Bolzano's intuition and deep insight seem to have preceded his ability to formulate the theory of measurable numbers with complete precision. 
The great debate as to whether measurable numbers is indeed a construction of real numbers has not yet been convincingly concluded. An affirmative answer would mean that Bolzano's construction pre-dates Cantor's and Dedekind's constructions by about forty years.  
Scholars' opinions vary; ranging from condemning the theory as inconsistent \citep{vanrootselaarBolzanoTheoryReal1964}, to classifying it as a precursor to Cantor's theory of real numbers \citep{rychlik1957Theorie}, \citep{sebestikLogiqueMathematiqueChez1992}, \citep{rusnock2000Bolzanos}, 
to believing that one small change in the initial definition can completely correct the entire theory \citep{laugwitz1965Bemerkungen}, \citep{russ2016Bolzanosa}, or as an attempt of a rigorous treatment of the real numbers still using some traditional notions such as variable numbers \citep{bellomoSumsNumbersInfinity2021}, a not-yet-modern conception of mathematics and numbers \citep{fuentesguillen2022Bolzanosa}.   

However, I do not want to continue this discussion now, but to point on the connection between Bolzano's measurable numbers and his concept of continuum. Therefore, I will first briefly outline the construction.

\subsection{Infinite number expressions}\label{b}

Bolzano starts with \emph{infinite number expressions} which are formed from rational numbers by using infinitely many arithmetic operations. He gives some examples
$$A = 1 + 2 + 3 + 4 + \dots \mbox{in inf.}$$
$$B = \frac{1}{2} - \frac{1}{4} + \frac{1}{8} - \frac{1}{16} + \dots \mbox{in inf.}$$
$$C = (1 -  \frac{1}{2})(1 - \frac{1}{4})(1 - \frac{1}{8})(1 - \frac{1}{16}) \dots \mbox{in inf.}$$
$$D = a  +   \frac{b}{1+1+1+\dots \mbox{in inf.}},$$ where $a, b$ is a pair of whole numbers. 

Bolzano distinguishes several kinds of infinite number expressions.  
\begin{itemize}
\item \emph{Measurable numbers} are the expressions which \enquote{we determine by approximation or the measurement as far as we please}. Formally, it means that for each natural number $q$ there is an integer $p$ and two positive number expressions $P_1, P_2$, the former possibly denoting zero, such that 
$$S = \frac{p}{q} + P_1 = \frac{p+1}{q} - P_2.$$
For instance, the expressions $C$, $D$ are measurable.  
\item \emph{Infinitely small positive} number are such that for each $q$ there are two positive number expressions $P_1, P_2$, the former possibly denoting zero, such that 
$$S =  P_1 = \frac{1}{q} - P_2.$$ 

For example $C$. \emph{Infinitely small negative} numbers are defined analogically. 
\item \emph{Infinitely great positive} numbers are greater than any $n$, for instance $A$. 
\end{itemize}

Bolzano proves many theorems on these kinds of numbers. Rational numbers are measurable. Measurable numbers are closed on addition and multiplication. Infinitely small numbers are closed on multiplication and addition. Infinitely great numbers are also closed on multiplication. The sum of finitely many infinitely small numbers is infinitely small. The product of an infinitely small (resp. great) and measurable number is infinitely small (resp. great), etc.

\subsection{Equality (equivalence)}  

The last step of the construction is the introduction of the \emph{equality}. Bolzano already explained his concept of \emph{equality} in BG II, \S 1 and \S 2. Two things are \emph{equal} if and only if their determining pieces are equal. They can be different which means they are not \emph{identical}, simply not the same. 
As far as measurable numbers are concerned Bolzano sometimes uses the term \emph{equivalence}, instead of \emph{equality}, which is more intelligible to our ears. 

The idea is that two measurable numbers $A,B$ are \emph{equal (equivalent)} iff they \enquote{have the same characteristic in the process of measuring}, i.e. they have the same measuring fractions. This original definition brought a problem as Bolzano himself probably recognized and changed it. According to the corrected definition  $A,B$ are \emph{equal} if and only if 
their difference $|A-B|$ \enquote{has the same characteristic as zero itself in the process of measuring}
$$A = B \textrm{\ iff \ }  |A-B| \textrm{ is infinitely small.}$$ 

Consequently, all infinitely small numbers are equal to $0$. If $A$ is measurable and $J$ is infinitely small then $A + J = A$. 
The ordering of measurable numbers is introduced on the same principle. $A$ is \emph{greater} than $B$,
$$A > B \textrm{\ iff \ } A - B \textrm{ is positive and not infinitely small.}$$

\subsection{Properties of measurable numbers}

After introducing equality and ordering, Bolzano has got \enquote{new} measurable numbers.  However, he does not distinguish them from the original ones and deals with their properties.   He proves that the ordering is transitive [\S 63], compatible with addition [\S 67], unbounded [\S 70],  linear [\S 73], Archimedean [\S 74], and dense [\S 79]. The multiplication is commutative and associative [\S 99], and the distributive law [\S 101] applies.  

So, measurable numbers have almost all important properties of contemporary real numbers. The only left is the completeness. This is what the three main theorems [\S\S 107, 109, 110] are for, see Section \ref{main}. But first we show how Bolzano's measurable numbers can be interpreted.

\section{The interpretation}\label{interpretation}

Lakatos in his paper about the significance of Cauchy for the analysis of continuum recalls Bolzano. 
\begin{quote}It is a most interesting historical fact that Bolzano, the best logical mind of the generation, made a real effort to clear the matters. He was possibly the only one to see the problems related the difference between the two continua: the rich Leibnizian continuum and, as he called it, its 'measurable' subset - the set of Weierstrassian real numbers. Bolzano makes it very clear that the field of 'measurable numbers' constitutes only an Archimedean subset of a continuum enriched by non-measurable - infinitely small and infinitely large quantities... No doubt, since Robinson has shed new light on the latter, historians will approach the Bolzano manuscript with new eyes and the relation between Bolzano's measurable numbers and non-measurable quantities and Robinson's standard and non-standard numbers will be clarified. \citep[p. 154]{lakatos1978Cauchy}. 
\end{quote}

That is the question I would like to answer, what is the relationship among Bolzano's measurable numbers, Cantor's real numbers and Robinson's non-standard numbers. 

It is most appropriate to interpret Bolzano's infinite expressions as sequences of partial results  arithmetic operations and indeed, most interpreters do. We get sequences of rational numbers. Measurable numbers correspond to Bolzano-Cauchy sequences (BC-sequences). Infinitely small numbers to sequences converging to $0$, and infinitely great numbers to divergent sequences. With a small correction in the definition of measurable numbers,\footnote{The problem is Bolzano's own example 
$$B = \frac{1}{2} - \frac{1}{4} + \frac{1}{8} - \frac{1}{16} + \dots \mbox{in inf.}$$ measurable? See \ref{b}. Its interpretation is the sequence of partial results oscillating around $\frac{1}{3}$, therefore for $q=3n$ the expression $B$ is neither an element of $[\frac{n}{3n}, \frac{n+1}{3n})$ nor $[\frac{n-1}{3n}, \frac{n}{3n})$, where $n \in \mathbb N$. This leads to a number of problems based on the uncertainty of whether or not the expressions leading to oscillating sequences are measurable. If not, then the measurable numbers are not closed to the sum. If the answer is yes, then Bolzano's definition of measurability is imprecise. It turns out, however, as already \citet[p. 407]{laugwitz1965Bemerkungen} suggested, that it is sufficient to change a little the definition of measurable numbers so that an infinite number expression $S$ is \emph{measurable} if for each natural number $q$ there is an integer $p$ and two positive number expressions $P_1, P_2$ such that 
$$S = \frac{p-1}{q} + P_1 = \frac{p+1}{q} - P_2.$$

This solution is also indicated by Bolzano's claims that an infinitely small measurable number $J$ equals $0$ and that both measurable numbers and infinitely small quantities are closed on addition.}
\emph{cum grano salis}, it is a consistent structure.

\subsection{The Cantor theory}

Cantor works with BC-sequences that he calls \emph{fundamental}. With each such sequence ${(a_n)}$ he associated a definite symbol $a$ as its limit. He then defined an equality of two fundamental sequences so that their difference is the sequence converging to zero. Similarly the ordering and the basic arithmetic operations. Finally, Cantor has got the structure of real numbers. 
\begin{quote} The real numbers take the form of purely mental symbols associated to the fundamental sequances. \citep[p. 128]{ferreirosjoseLabyrinthThought2007}\end{quote} 

\begin{center} 
\begin{tabular}{|c|c|}
 \hline
 \bf The Bolzano theory & \bf The Cantor theory \\
\hline
infinite number expressions &  sequences of rational numbers\\

infinitely small numbers  &  sequences converging to zero \\
infinitely large numbers  &  divergent sequences \\
measurable numbers  &  BC-sequences
\\$x,y$ are measurable numbers &  $(x_n)$, $(y_n)$ are BC-sequences \\

$x - y$ is infinitely small  & $\lim_{n\to\infty}(x_n - y_n) = 0$ \\
$x = y$ & $(x_n) = (y_n)$  \\
\enquote{new} measurable numbers & real numbers \\
\hline
 \end{tabular}
 \end{center} 

Formally, it corresponds. Nevertheless, there is a conceptual difference between Cantor's fundamental sequences and Bolzano's measurable numbers.  Bolzano himself called and considered them \emph{numbers} and \emph{quantities}. He proved many propositions on measurable numbers. Together with infinitely small and infinitely large \emph{numbers} they form a rich non-Archimedean continuum of the Leibnizian type. However, in Cantor's conception there is no place for it.  

\subsection{The non-standard theory}

Besides the well-known Cantor and Dedekind construction of real numbers, there is a construction using the non-standard model of \emph{rational} numbers \citep[p. 14]{albeverio2009Nonstandardb}. It works with sequences of rational numbers $\mathbb Q^\mathbb N$ and with a non-principal ultrafilter $\mathcal U$ on $\mathbb N$. The ultraproduct $\mathbb Q^* = \mathbb Q^\mathbb N/\mathcal U$ contains infinitely small numbers $\mathbb Q_i$ and infinitely great numbers defined as usually. Finite numbers $\mathbb Q_f$ are those that are not infinitely great. 
$$Q^\mathbb N/\mathcal U =  \mathbb Q^* \supseteq \mathbb Q_f \supseteq \mathbb Q_i .$$

Two finite non-standard numbers $x,y \in \mathbb Q_f$ are infinitely close $x \approx y$ if their difference is infinitely small. By the factorization modulo infinite closeness $\approx$ we obtain the complete linearly ordered Archimedean field of real numbers $$\mathbb R = \mathbb Q_f/\approx.$$

Real numbers are represented by factor-classes called \emph{monads}. So a monad is a class of infinitely close non-standard rational numbers.

Both finite non-standard rational numbers $\mathbb Q_f$ and Bolzano's \enquote{old} measurable numbers form a non-Archimedean commutative ring. Both non-standard infinitely small rational numbers $\mathbb Q_i$ and Bolzano's infinitely small numbers form their maximal ideals. The infinite closeness $\approx$ is defined exactly as Bolzano's equality, it is an equivalence relation in both cases. 

According to the well-known algebraic theorem, the factorization of a commutative ring modulo its maximal ideal yields a field, the same field as we get by the factorization modulo $\approx$. In case of $\mathbb Q_f/\approx$, it is a field of real numbers. Bolzano proved that his field is dense, linearly ordered, complete, Archimedean, and contains rational numbers. This uniquely determined structure is isomorphic to real numbers.

Bolzano does not distinguish between \enquote{old} measurable numbers and their factor-classes, \enquote{new} measurable numbers. It is a correct, albeit unusual, way. Arithmetic properties of \enquote{old} measurable numbers are preserved while \emph{equality} and \emph{order} are defined only for their factor-classes.

\begin{center} 
\begin{tabular}{|c|c|}
 \hline
 \bf The Bolzano theory & \bf The non-standard theory \\
\hline
infinite number expressions &  sequences of rational numbers $\mathbb Q^\mathbb N$ \\ measurable numbers   &  finite rational numbers $\mathbb Q_f$ \\
infinitely small numbers  &  infinitely small numbers $\mathbb Q_i$  \\
infinitely great numbers  &  infinitely great numbers \\
the equality $x = y$  & the equivalence $x \approx y$ \\
measurable numbers with $=$ & $\mathbb Q_f/\approx \ = \ \mathbb Q_f/\mathbb Q_i \ = \  \mathbb R$. \\
\hline
 \end{tabular}
 \end{center} 

Of course, this comparison is not accurate. The ultraproduct $\mathbb Q^*$ is constructed with help of an ultrafilter $\mathcal U$ which is a non-constructive, \emph{intangible} object, we need the \emph{Axiom of Choice} to prove its existence. However, this makes the ultraproduct $\mathbb Q^*$ as well as $\mathbb Q_f$ linearly ordered which cannot be said of Bolzano's \enquote{old} measurable numbers. On the other side, measurable numbers are built constructively and do not require any intangible object. However, the result is the same, real numbers and measurable numbers are isomorphic.

This interpretation demonstrates two things.  First, borrowing the notion of a monad from non-standard theory, each \enquote{new} measurable number is not a mere symbol, it is a \emph{monad} of equal \enquote{old} measurable numbers. Two measurable numbers are \emph{equal} if they \enquote{have the same properties in the process of measuring}. 
$$0 = \frac{1}{1+1+1+\dots \mbox{in inf.}} = \frac{1}{3+3+3+\dots \mbox{in inf.}} = \frac{1}{1+2+3+\dots \mbox{in inf.}}$$ 

Second, by comparison with Robinson's non-standard analysis, we show why Bolzano's infinitely small numbers are not suitable for the infinitesimal calculus.

%

\subsection{The Robinson non-standard analysis}

The Robinson theory uses a similar construction as above. The only difference is that it works with sequences of \emph{real} numbers. However, its goal is not the arithmetization of continuum but the reconstruction of infinitesimal calculus in a non-standard analysis. 

Due to the ultrafilter $\mathcal U$ the ultraproduct $\mathbb R^* = \mathbb R^\mathbb N/\mathcal U$ is a linearly ordered non-Archimedean field that is an elementary extension of $\mathbb R$ where the \emph{Transfer Principle} holds. It means that the same properties of the first-order logic hold for elements of $\mathbb R$ and $\mathbb R^*$. It enables to define the same structures in $\mathbb R$ and in $\mathbb R^*$ and guarantees their correspondence. Thus, functions defined on $\mathbb R$ can be naturally extended on $\mathbb R^*$. Limits, continuity, derivative and integrals can be defined in a simple and intuitive way based of infinitesimals.



However, the \emph{Transfer Principle} does not hold between $\mathbb R$ and $\mathbb Q^*$. Functions defined on real numbers cannot be simply extended on non-standard rational numbers. Therefore although $\mathbb Q^*$ contains infinitesimals the structure $\mathbb Q^*$ is not suitable for non-standard analysis.
 
The same argument holds for Bolzano's measurable numbers. Functions defined on his \enquote{new} measurable numbers cannot be simply extended on his \enquote{old} measurable numbers containing infinitely small quantities. To define an infinitesimal calculus on these quantities would be very difficult. 

Bolzano was probably aware of such difficulties. In the following \emph{Theory of Functions} \citep[pp. 429 - 589]{russ2004Mathematical} he returned to his primordial idea of building the calculus on the basis of \enquote{quantities which can become smaller than any given quantities} which is the idea of the later Weierstassian $\epsilon - \delta$ calculus. Nevertheless, he did not reject infinitely small and infinitely large quantities entirely. In \emph{Paradoxes of the Infinite} he defended their existence as well as the existence of infinitely great quantities as meaningful notions.

\section{Completeness}\label{main}

Bolzano's three main theorems in [\S\S 107, 109, 110] are the culmination of RZ, \enquote{les trois grands théor\`emes par lesquels culmine la théorie de Bolzano} \citep[p. 384]{sebestikLogiqueMathematiqueChez1992}. They prove the three forms of completeness of real numbers.
\begin{itemize} 
\item \emph{Cauchy completeness.} Every BC-sequence of real numbers has a limit. 
\item \emph{Supremum completeness.} Every non-empty subset of real numbers having an upper bound has a supremum.
\item \emph{Dedekind completeness.} Every Dedekind cut of rational numbers is generated by a certain real number. 

\end{itemize}

Moreover, I believe that the main purpose of the last theorem [\S 110] is to prove that measurable numbers satisfy the characteristic property of the continuum.

\begin{itemize}
\item \emph{Bolzano completeness.} Every measurable number has at least one neighbour at every distance sufficiently small.
\end{itemize}

\subsection{The Bolzano-Cauchy Theorem - [RZ \S 107]}

The theorem [RZ \S 107] is clearly the Bolzano-Cauchy.  
While slightly different wording, it has the same meaning as  [RB \S 7]: every BC-sequence of measurable numbers has a measurable limit. In 1817, the existence of a limit could not be proved  because of the non-existence of the arithmetic description of continuum. The second proof in RZ is not yet perfect, lacking some details, but is nevertheless almost successful. \citep[p. 188]{rusnock2000Bolzanos}, \citep[p. 94]{bellomoSumsNumbersInfinity2021}. This theorem is the key to proving the following more general results.

\subsection{The Supremum Theorem - RZ \S 109}\label{supremum}

This theorem is variously called. According to \citet[p. 103]{sebestikLogiqueMathematiqueChez1992}, \enquote{la théoreme qui affirme l'existence d'une \emph{borne supérieure}}; \citet[p. 150]{russ2004Mathematical} says that \enquote{it could be called  the \emph{greatest lower bound property}}; and by  \cite[p. 210]{rusnock2020Bolzano},, \enquote{it is closely related to the \emph{least upper bound property}}. 
Again, almost the same theorem is already stated and proved in [RB \S 12].  All commentators agree that its proof is well organized, clear and rigorous, though long. 
I take the liberty of calling it the \emph{Supremum Theorem} and will show that this is precisely its meaning.


\begin{quote} If we know about a certain property $\mathcal B$, that it belongs, not to all values of a variable measurable number $X$ which are greater (or smaller) than a certain
value $U$, but to all which are smaller (greater) than $U$, then we can definitely claim
that there is a measurable number $A$ which is the greatest (smallest) of those of
which it can be said that all smaller (greater) $X$ have the property $\mathcal B$. It is left still undecided here whether the value $X = A$ itself also has this property. [RZ \S 109].
\end{quote}

This theorem is unusually worded. I reformulate it for better comprehension. 

\begin{enumerate}
\item In symbolic notation (the case outside parentheses): 

 Let $\mathcal B(X)$ means that \enquote{$X$ has the property $\mathcal B$}. Let $U$ be a measurable number such that
$$(\ast) \quad (\forall X)(X < U \Rightarrow \mathcal B(X)) \wedge (\exists X)(X > U \wedge \neg \mathcal B(X).$$ 

Then there is a measurable number $A$ such that 
$$(\ast \ast) \quad (\forall X)(X < A \Rightarrow \mathcal B(X) \wedge (\forall Y)((\forall X)(X < Y \Rightarrow \mathcal B(X)) \Rightarrow Y \leq A).$$

\item We use the law of \emph{contraposition} and the \emph{linearity}\footnote{From linearity of measurable numbers which Bolzano proved in RZ \S 73 follows that $\neg(X > Y)$ iff $X \leq Y$.} of the ordering:  

Let $\mathcal B(X)$ means that \enquote{$X$ has the property $\mathcal B$}. Let $U$ be a measurable number such that  

$$(\ast) \quad  (\forall X)(\neg \mathcal B(X) \Rightarrow X \geq U) \wedge (\exists X)(X > U \wedge \neg \mathcal B(X) )$$

which means that $U$ is a lower bound of the non-empty set $\{x; \neg \mathcal B(x)\}$.
Then there is a measurable number $A$ such that
$$(\ast \ast) \quad (\forall X)(\neg \mathcal B(X) \Rightarrow X \geq A) \wedge (\forall Y)((\forall X)(\neg \mathcal B(X) \Rightarrow X \geq Y) \Rightarrow Y \leq A)$$ 

which means that $A$ is the greatest lower bound, the \emph{infimum} of $\{x; \neg \mathcal B(x)\}$.
$$A = \inf\{X; \neg \mathcal B(X)\}.$$ 

The theorem [\S 109] states that if the non-empty set $\{x; \neg \mathcal B(x)\}$ has a lower bound then it has an infimum, exactly as \citet[p. 269]{russ2004Mathematical} comments.

\item 
Let us denote the measurable number $U$ for which $(\forall X)(X < U \Rightarrow \mathcal B(X))$, i.e. a lower bound of $\{X; \neg \mathcal B(X)\}$, as the upper \emph{boundary} of $\mathcal B(X)$. Then \S 109 states:

 If there is an upper (or lower) boundary of $\mathcal B(X)$ and $\mathcal B(X)$ does not apply for all measurable numbers then there is the greatest upper (smallest lower) boundary $A$. It is left undecided whether $\mathcal B(A)$ or not. 

\end{enumerate}

The last formulation also corresponds to the way of Bolzano's proof. He constructs a sequence of still greater upper boundaries of $\mathcal B$ and demonstrates that either some of them is the greatest one or that this sequence has a BC property. Hence, according to [\S 107] it has a limit that is the greatest upper boundary of $\mathcal B$. 

\subsection{The Intermediate Value Theorem - [RB \S 15]}\label{IMV}

Let us return briefly to the \emph{Purely Analytic Proof} (RB) from 1817. The \emph{Bolzano-Cauchy} [\S 7] and the \emph{Supremum} [\S 12] Theorems are followed by  the \emph{Intermediate Value Theorem} [\S 15], which is the principal result of this essay. 

\begin{quote} If two functions $f(x)$ and $\varphi(x)$, vary according to the law of continuity \dots
and furthermore if $f(\alpha) < \varphi(\alpha)$ and $f(\beta) > \varphi(\beta)$, then there is always a certain value of $x$ between $\alpha$ and $\beta$ for which $f(x) = \varphi(x)$. [RB \S 15]. \end{quote}

Bolzano uses repeatedly the corollary of his definition of continuity of a function, see Section \ref{continuous}, saying that if functions $f(x), \varphi(x)$ are continuous then $$(\ast) \quad f(\alpha) < \varphi(\alpha) \Rightarrow f(\alpha + i) < \varphi(\alpha + i)$$ 

provided $i$ is taken sufficiently small. 
\medskip

In the proof of  [RB \S 15], Bolzano defines the property $\mathcal B$  
$$\mathcal B(i): f(\alpha + i) < \varphi(\alpha + i).$$ 

According to $(\ast)$, $\mathcal B(i)$ 
has an upper boundary. Nevertheless, it does not hold for all values, it does not hold for $i = \beta - \alpha$. Therefore according to [RB \S 12], $\mathcal B(i)$ has the greatest upper boundary $u$.  

 
For this value $u$ itself, $\mathcal B(u)$, i.e. $f(\alpha + u) < \varphi(\alpha + u)$, cannot be valid, otherwise by $(\ast)$ there would be $\omega$ such that $\mathcal B(u + \omega)$, so $u$ would not be the \emph{greatest} upper boundary. But still less can be that $f(\alpha + u) > \varphi(\alpha + u)$ because then also $f(\alpha + u - \omega) > \varphi(\alpha + u - \omega)$, hence $u$ is not an \emph{upper boundary}. Consequently, $$f(\alpha + u) = \varphi(\alpha + u).$$ 

\subsection{The Bolzano Completeness - [RZ \S 110]}\label{BC}

Now, let us ask ourselves whether and why [RZ \S 110] is the culmination of the previous theorems as [RB \S 15]. According to \citet[p. 385]{sebestikLogiqueMathematiqueChez1992}, it is resembling Dedekind's cuts, but unlike them it is not a \emph{construction} of real numbers, but a demonstration of a property of real numbers \emph{already constructed}.

\begin{quote} If the variable but measurable number Y always remains greater than
the variable but measurable number X, and if also there is no greatest value of the former and no smallest value of the latter, then there is always at least
one measurable number A which always lies between the two limits X and Y.
Furthermore, if the difference Y - X cannot decrease indefinitely, then there are
infinitely many such measurable numbers lying between X and Y. But if this
difference does decrease indefinitely, then there is only a single [number]. Finally
if the difference Y - X decreases indefinitely and either X has a greatest value,
or Y has a smallest value, then there is not a single measurable number which
always lies between X and Y. [RZ \S 110].
\end{quote}

The notion of a \emph{variable} and a \emph{constant} quantity was commonly used, for instance Guillaume de l'Hospital defines: 
\enquote{Those quantities are called variable which increase or decrease
continually, as opposed to constant quantities that remain the same while others
change.} \citep[p. 139]{shapiroHistoryContinuaPhilosophical2021}. Thus, a variable measurable number can take any value on a set defined by some property. We can reformulate [RZ \S 110] in the contemporary mathematical language. 

\begin{quote} 
Let $X, Y$ be two subsets of measurable numbers such that 
$$(\forall x \in X)(\forall y \in Y)(x < y).$$
Moreover, $X$ does not have a greatest element and $Y$ does not have a smallest element. 

Then there is a measurable number $a$ which \emph{is between $X$ and $Y$} $$(\forall x \in X)(\forall y \in Y)(x < a < y).$$ 
\begin{enumerate}[(1)]
\item If $(\exists n \in \mathbb N)(\forall x \in X)(\forall y \in Y)(y - x > \frac{1}{n})$ then there are infinitely many such numbers. 
\item  If $(\forall n \in \mathbb N)(\exists x \in X)(\exists y \in Y)(y - x < \frac{1}{n})$ then there is exactly one such number.
\item If moreover either $X$ does have the greatest element or $Y$ does have the smallest element then there is no such number. 
\end{enumerate}
\end{quote}

The beginning of Bolzano's proof is little complicated because of his term of \emph{boundary}. He considers any $u \in X$ as an upper boundary of the complement of $Y$.
$$(\forall x)(x < u \Rightarrow x \notin Y) \wedge (\exists x)(x > u \wedge x \in Y),$$
which actually is a lower bound of $Y$. 
$$(\forall x)(x \in Y \Rightarrow x \geq u) \wedge (\exists x)(x \in Y \wedge  x > u).$$
According to the \emph{Supremum Theorem} there is $a$, the greatest upper boundary of the complement of $Y$, i.e. the greatest lower bound, infimum, of $Y$ which has the desired property of being between $X$ and $Y$. 
$$a = \inf\{x; x \in Y\}.$$
In case of (1), Bolzano demonstrates that also $a-\frac{1}{n}$ lies between $X$ and $Y$. 
From the density of measurable already proven in [RZ \S 84] it follows there are infintely many measurable numbers between $a - \frac{1}{n}$ and $a$. Bolzano proves (2) by contradiction. If (3) then the greatest element of $X$ or the smallest element of $Y$ must be equal to $a$. There is no other number lying between $X$ and $Y$.

Of course this theorem also implies the \emph{Dedekind completeness}. A Dedekind cut corresponds either to the case (2) or (3). In both cases, there is exactly one measurable number that generates it.

Nevertheless, I believe that the most important corollary of [\S 110] is the \emph{Bolzano completeness} of measurable numbers.   

\begin{quote} Measurable numbers are so arranged that every single of them has at least one neighbour at every distance however small.
 \end{quote}

\emph{Proof.} Let $X$ be a variable measurable number less than a given distance $d$, $Y$ be a variable number greater than $d$. Surely $X < Y$ and the difference $Y - X$ does decrease indefinitely. By (2), there is exactly one measurable number $a$, lying between $X$ and $Y$. Neither $a$ is less nor greater than the distance $d$. The only possibility is that $a$ equals the given distance. There is an option that a variable measurable number $X$ has a greatest element, for instance when we consider $X$ less \emph{or equal} than $d$. Then according to (3), the greatest element is equal to $d$ and there is no other such number.

\section{Back again}\label{back}

From the very beginning, Bolzano considered points of a line as abstract objects. 

\begin{quote}
The concept of point, as a mere \emph{characterization of space that is itself no part of space - cannot be dispensed in with} geometry. The point is a merely imaginary object. (BG II \S 5).
\end{quote}

Measurable numbers as imaginary objects characterize not only the points of a line, but also places (points) where the substances of which any spatial or temporal continuum is composed occur.

\begin{quote} For me space, in a similar way as time, is \emph{not an attribute} of substances but only a \emph{determination} of them. Indeed, I call those determinations of created substances \dots the \emph{places} at which they occur. The collection of all places I call \emph{space}, the whole of space. [PU \S 40]. \end{quote}

Only now when Bolzano can characterize continuous extensions by measurable numbers, he can confidently clarify the properties of continuum and thus reject the paradoxes associated with them.
\begin{quote} In \emph{Paradoxes}, at the end of his life, it is no longer even the continuity of the real line that is analysed in principle, but a notion of continuity applicable to any space in general. \citep[p. 521]{granger1969Concept}.\footnote{\enquote{Dans les \emph{Paradoxes}, a la fin de sa vie, ce n'est plus même la continuité de la droite réelle qui est en principe analysée, mais une notion de continuité applicable a toute espace d'ensemble en général.}}  \end{quote}

Bolzano's mature concept of continuum\footnote{In later writings, Bolzano uses the terms \emph{extension},  \emph{continuous extension} and occasionally \emph{continuum} as nearly equivalent, which was not always the case  \cite[p. 200]{rusnock2020Bolzano}.} is found in \emph{Paradoxes of the Infinite} [PU] written in 1847 - 48 and published posthumously in 1851, particularly in [\S 38]. His aim was to summarize general properties of the temporal, the spatial and even of the material continua as he writes at the beginning: \enquote{Therefore we shall consider all them together.} When we follow his arguments we see that they apply as well for measurable numbers.

\begin{enumerate}

\item As before, he emphasises that any continuous extension is composed of simple parts that lack extension.
\begin{quote}
Continuous extensions are composed of parts that have no extension but are absolutely simple (points in time or space, atoms, i.e. simple substances in the universe within the realm of reality).
\end{quote}
An interval, that has an extension, is composed from measurable numbers, that have no extension. 

\item This infinite multitude of simple parts must be arranged in a manner that is now called  \emph{density}: between any two of them, there exists an infinite multitude of others.

\begin{quote} That every two instants in time are separated by an infinite multitude of instants in between, that likewise between every two points in space there is an infinite multitude of points lying in between them, that even in the realm of reality between every two substances there is an infinite multitude of others - is of course conceded. \end{quote}

\enquote{If $A$ and $C$ are a pair of unequal measurable number then there is always a third measurable number which lie between the two.}  [RZ \S 79].

\item Only an infinite multitude of parts can give rise to an extension.  
\begin{quote} This much follows, that with two simple parts alone, or even with three, four or any merely finite multitude of them, no extension is produced.  \end{quote}

\enquote{There is an infinite multitude of measurable numbers between any two different measurable numbers.} [RZ \S 84].

However, already these conditions are not sufficient to constitute a continuum. For instance, they hold for rational numbers, which does not provide a complete arithmetization of the line.\footnote{Bolzano also gives an example that the density is not a sufficient condition.  
\enquote{For consider two points $a,b$ along with a point $c$ in the middle between them, and points $d,e$ in the middle, respectively, between $a, c$ and $c, b$, and so on ad infinitum. None of the points belonging to this system, he notes, has a closest neighbour; yet no geometer would claim that these points are connected with each other and form an extension. For even though each point has neighbours that are arbitrarily close to it, for each there are also arbitrarily small distances at which it has no neighbour. }[BBGA 2A.11/1, \S6], cited according to  \cite[p. 201]{rusnock2020Bolzano}.}

\item Finally, the key property of continuity, the \emph{Bolzano completeness}. 
\begin{quote}
A continuum exists where, and only where, a collection of simple objects (of points in time or space or even of substances) occurs which are so arranged that every single one
of them has at least one neighbour in this collection at every distance however
small.
\end{quote}

It follows from [RZ \S 110] as we proved in Section \ref{BC}. 

The delicacy of this definition is well illustrated by the example in [PU \S 41, 3]. Let us consider a straight line $az$, let $a_1$ be its midpoint, $a_2$ be the midpoint of $a_1z$,  $a_3$ the midpoint of $a_2z$, and so on. When we omit all points $a_1, a_2, \dots$ and the point $z$ then the collection of remaining points still deserves the name a \emph{line}. But when include $z$ in this collection then it is not a continuous extension. The point $z$ has no neighbours in all distances $\frac{az}{2^n}$.

\end{enumerate}

Bolzano's topological concept of the continuum is well illustrated by the following excerpt which is just as appropriate for his measurable numbers, especially as interpreted in the non-standard model.   

\begin{quote}  
Points are simple parts of space, they therefore have no boundaries, no
right and left sides, no extension. If one had only a part in common with the other then it would be absolutely the same as it, and if it is to have something different from it, then both must lie completely outside one another and there must therefore be space
for another point lying between them. Indeed because the same holds of these intermediate points in comparison with those two, [there is space] for an infinite multitude of points. [PU \S 38].
\end{quote}

\section{Conclusion}

Bolzano's conception of the continuum was corpuscular from beginning to end. He established a clear condition, the \emph{Bolzano Completeness}, under which a collection of points forms a continuum: no point is isolated, i.e. it has at least a sufficiently small neighborhood containing no gaps. The introduction of measurable numbers meant the necessary arithmetization of continuum. It was important to show they form a continuum, that the \emph{Bolzano Completeness} applies.  

Measurable numbers interpreted in a non-standard model 

If we were going to argue about what Bolzano is ahead in, definition of continuity of a function, BC-criterion, the arithmetization of continuum, completness theorems. 

\bibliography{Bibliography}

\end{document}